\documentclass[11pt, a4paper]{article}
\usepackage{amsfonts}
\usepackage{mathrsfs}
\usepackage{amssymb}
\usepackage{amsmath}
\begin{document}

\title{The strong converse inequality for de la Vall\'{e}e Poussin means on the sphere
\thanks{The research was supported by the National
Natural Science Foundation of China (No. 60873206)}}

\author{Ruyue Yang  \and  Feilong Cao\thanks{Corresponding author: Feilong Cao,
E-mail: \tt feilongcao@gmail.com}
 \and Jingyi Xiong}

  \date{}
\maketitle

\begin{center}
\footnotesize
  Department of Mathematics, China
 Jiliang University,

 Hangzhou 310018, Zhejiang Province,  P R China

\begin{abstract}
This paper discusses the approximation by de la Vall\'{e}e Poussin means $V_nf$ on the unit sphere. Especially, the lower bound of approximation is studied. As a main result, the strong converse inequality for the means is established. Namely, it is proved that  there are constants $C_1$ and $C_2$ such that
\begin{eqnarray*}
C_1\omega\left(f,\frac{1}{\sqrt n}\right)_p \leq \|V_{n}f-f\|_p \leq C_2\omega\left(f,\frac{1}{\sqrt n}\right)_p
\end{eqnarray*}
for any  $p$-th Lebesgue integrable or continuous function $f$ defined on the sphere, where $\omega(f,t)_p$ is the modulus of smoothness of $f$.
\end{abstract}
\end{center}

\section{Introduction}

Motivated by geoscience,  meteorology and oceanography,  sphere-oriented mathematics has gained increasing attention in recent decades. As the main tools, spherical positive polynomial operators play prominent roles in the approximation and the interpolation on the sphere by means of orthonormal spherical harmonics.
Several authors such as Ditzian \cite{Ditzian2004}, Dai and Ditzian \cite{Dai2009}, Bernes and Li \cite{Bernes1968}, Wang and Li \cite{Wang2000}, Nikol'ski\v{\i} and Lizorkin \cite{Nikolskii1985, Lizorkin1983} introduced and studied some spherical versions of some known one-dimensional polynomial operators, for example, spherical Jackson operators (see \cite{Lizorkin1983}), spherical de la Vall\'{e}e Poussin operators (see \cite{Berens1993}, \cite{Wang2000}), spherical delay mean operators (see \cite{Stein1957}) and best approximation operators (see \cite{Ditzian2004}, \cite{Dai2009}, \cite{Wang2000}) etc..

The main aim of the present paper is to study the  approximation by the de la Vall\'{e}e Poussin means on the unit sphere.

For to formulate our results, we first give some notations. Let $\mathbb R^{d}$, $d\ge 3$, be the Euclidean space of the points $x:=(x_1,x_2,\dots, x_d)$ endowed with the scalar product
$$
x\cdot x^{\prime} = \sum_{j=1}^d x_jx_j^{\prime}\quad (x, x^{\prime}\in \mathbb{R}^{d})
$$
and let $\sigma:=\sigma^{d-1}$ be the unit sphere in $\mathbb R^d$ consisting of the points $x$ satisfying $x^2 =x\cdot x=1$.

We shall denote the points of $\sigma$ by $\mu$, and the elementary surface piece on $\sigma$ by $d\sigma$. If it is necessary, we shall write $d\sigma:\equiv d\sigma(\mu)$ referring to the variable of integration. The surface area of $\sigma^{d-1}$ is denoted by $|\sigma^{d-1}|$, and it is easy to deduce that
$$
|\sigma^{d-1}| = \int_{\sigma}d\sigma = \frac{2\pi^{\frac{d}{2}}}{\Gamma(\frac{d}{2})}.
$$

By $C(\sigma)$ and $L^p(\sigma)$, $1\leq p< +\infty$, we denote the space of continuous, real-value functions and the space of (the equivalence classes of ) $p$-integrable functions defined on $\sigma$ endowed with the respective norms
$$
\|f\|_{\infty} := \max_{\mu \in \sigma }|f(\mu)|,
$$

$$
\|f\|_{p} := \left(\int_{\sigma}|f(\mu)|^{p}d\sigma(\mu)\right)^{1/p},\quad 1 \leq p < \infty.
$$
In the following, $L^p(\sigma)$ will always be one of the spaces $L^p(\sigma)$ for $1\leq p<\infty$, or $C(\sigma)$ for $p=\infty$.

Now we state some properties of spherical harmonics (see \cite{Wang2000}, \cite{Freeden1998}, \cite{Muller1966}). For
integer $k \geq 0$, the restriction of a homogeneous harmonic polynomial of degree $k$ on the unit sphere is
called a spherical harmonic of degree $k$. The class of all spherical harmonics of degree $k$ will be denoted by
$\mathcal {H}_{k}$, and the class of all spherical harmonics of degree $k \leq n$ will be denoted by $\Pi_{n}^{d}$. Of course,
$\Pi_n^{d} = \bigoplus_{k=0}^{n}\mathcal{H}_{k}$, and it comprises the restriction to $\sigma$  of all algebraic polynomials in $d$ variables of total degree not exceeding $n$. The dimension of $\mathcal{H}_{k}$ is given by
\[
d_k^{d} := \mbox{dim}\; \mathcal{H}^{d}_k := \left\{\begin{array}{ll}
\displaystyle\frac{2k+d-2}{k+d-2}{{k+d-2}\choose{k}}, & k\geq 1;\\
1, & k=0,
\end{array}
\right.
\]
and that of $\Pi_n^{d}$ is $\sum_{k=0}^{n} d^{d}_k$.

The spherical harmonics have an intrinsic characterization. To describe this, we first introduce the Laplace-Beltrami operators (see \cite{Muller1966}) to sufficiently smooth functions $f$ defined on $\sigma$, which is the restriction of Laplace operator
$\Delta:=\sum_{i=1}^{d}\frac{\partial^2}{\partial x_{i}^{2}}$ on the sphere $\sigma$, and  can be expressed  as
$$
Df(\mu) := \Delta f\left(\frac{\mu}{|\mu|}\right)\Big|_{\mu \in \sigma}.
$$
Clearly, the operator $D$ is an elliptic, (unbounded) self-adjoint operator on $L^2(\sigma)$, is invariant under arbitrary coordinate changes, and its spectrum comprises distinct eigenvalues $\lambda_k :=-k(k+d-2),\ k=0,1,\dots$, each having finite multiplicity. The space $\mathcal {H}_k$ can be characterized intrinsically as the eigenspace corresponding to $\lambda_{k}$, i.e.
$$
\mathcal{H}_{k} = \{\Psi \in C^{\infty}(\sigma): D\Psi = -k(k+d-2)\Psi \}.
$$
Since the $\lambda_{k}$'s are distinct, and the operator is self-adjoint, the spaces $\mathcal{H}_{k}$ are mutually orthogonal;
also, $L^2(\sigma) =$ closure $\{\bigoplus_k\mathcal{H}_k\}$. Hence, if we choose an orthogonal basis $\{Y_{k,l}: l = 1,\dots,d_{k}^{d}\}$ for each $\mathcal {H}_{k}$, then the set $\{Y_{k,l}: k=0, 1, \dots, l=1, \dots, d_k^{d}\}$ is an orthogonal basis for $L^2(\sigma)$.

The orthogonal projection $Y_k: L^1(\sigma)\rightarrow \mathcal{H}_{k}$ is given by
$$
Y_k(f;\mu) := \frac{\Gamma(\lambda)(k+\lambda)}{2\pi^{\lambda+1}} \int_{\sigma}P^{\lambda}_{k}(\mu\cdot\nu)f(\nu)d\sigma(\nu),
$$
where $2\lambda=d-2$, and $P^{\lambda}_k$ are the ultraspherical (or Gegenbauer) polynomials defined by the generating equation
$$
(1-2r\cos\theta +r^2)^{-\lambda} = \sum^{\infty}_{k=0}r^{k} P_{k}^{\lambda}(\cos\theta),\quad 0\leq\theta\leq \pi.
$$
The further details for the ultraspherical  polynomials can be found in \cite{Stein1971}.

For an arbitrary number $\theta$, $0<\theta<\pi$, we define the spherical translation operator of the function $f\in L^p(\sigma)$ with a step $\theta$ by the aid of the following equation (see \cite{Rudin1950}, \cite{Bernes1968}):
\begin{equation}\label{ch1eq1}
S_\theta(f) := S_{\theta}(f;\mu) := \frac{1}{|\sigma^{d-2}|\sin^{d-2}\theta}\int_{\mu\cdot \nu = \cos\theta} f(\nu) d\sigma(\nu),
\end{equation}
where $|\sigma^{d-2}|$ means the $(d-2)$-dimensional surface area of the unit sphere of $\mathbb R^{d-1}$. Here we integrate over the family of points $\nu\in\sigma$ whose spherical distance from the given point $\mu\in\sigma$ (i.e. the length of minor arc between $\mu$ and $\nu$ on the great circle passing through them) is equal to $\theta$. Thus $S_\theta(f;\mu)$ can be interpreted as the mean value of the function $f$ on the surface of  $(d-2)$-dimensional sphere with radius $\sin \theta$.

The properties of spherical translation operator (\ref{ch1eq1}) are well known; see e.g., \cite{Bernes1968}. In particular, it can be expressed as the following series
$$
S_{\theta}(f;\mu) = \sum_{k=0}^{\infty}\frac{P_{k}^{\lambda}(\cos \theta)}{P_{k}^{\lambda}(1)}Y_{k}(f;\mu) := \sum_{k=0}^{\infty}Q_{k}^{\lambda}(\cos\theta)Y_{k}(f;\mu)
$$
and for any $f\in L^p(\sigma)$
$$
\|S_\theta f\|_{p} \leq \|f\|_{p},
$$

$$
\lim_{\theta\rightarrow 0}\|S_\theta f-f\|_{p} = 0,
$$
where $Q_{k}^{\lambda}(\cos\theta) := \frac{P_{k}^{\lambda}(\cos \theta)}{P_{k}^{\lambda}(1)}$. We usually apply the translation operator to define spherical modulus of smoothness of a function $f \in L^{p}(\sigma)$, i.e. (see \cite{Wang2000})
$$
\omega(f,t)_p := \sup_{0< \theta \leq t}\|f-S_{\theta}f\|_{p}.
$$
Clearly, the modulus is meaningful to describe the approximation degree and the smoothness of functions on $\sigma$, which has been widely used in the study of approximation on sphere.

We also need a $K$-functional on sphere $\sigma$ defined by (see \cite{Ditzian2004}, \cite{Wang2000})
\begin{equation}\label{ch1eq2}
K(f; t)_p:=\inf\left\{\|f-g\|_{p}+t^{2}\|Dg\|_p: \ g,  Dg \in L^p(\sigma)\right\}. \quad 0<t<t_0.
\end{equation}
For the modulus of smoothness and $K$-functional, the following equivalent relationship has been proved (see \cite{Ditzian2004})
\begin{equation}\label{ch1eq3}
\omega(f,t)_p \approx K(f,t)_p.
\end{equation}
Here and in the following, $a\approx b$ means that there are positive constants $C_1$ and $C_2$ such that
$C_1\leq a\leq C_2b$. We denote by $C_i (i=1, 2, \dots)$ the positive constants independent of $f$ and $n$, and by $C(a)$ the positive constants depending only on $a$. Their value will be different at different occurrences, even within the same formula.

Define the kernel of de la Vall\'{e}e Poussin as
\begin{equation}\label{ch1eq4}
v_{n}(t)=\frac{1}{I_{n,d}}\Big(\cos \frac{t}{2}\Big)^{2n}, \quad n \in \mathbb N,
\end{equation}
where the constant $I_{n,d}$ satisfies
$$
\int_{\sigma}v_{n}(\stackrel\smile {\mu\nu})d\sigma(\nu)=|\sigma^{d-2}|,
$$
and $\stackrel\smile{\mu\nu}$ is the spherical distance between the points $\mu$ and $\nu$, i.e. the length of minor arc of great circle crossing $\mu$ and $\nu$. Then the convolution resulted by the kernel is
\begin{equation}\label{ch1eq5}
V_{n}(f;\mu)=(f*v_{n})(\mu)=\frac{1}{\big|\sigma^{d-2}\big|} ~\int_{\sigma}f(\nu)v_{n}(\stackrel\smile{\mu\nu}) d\sigma(\nu), \quad f\in L^1(\sigma),
\end{equation}
which is called de la Vall\'{e}e Poussin means on the sphere.

The means were introduced by de la Vall\'{e}e Poussin in 1908 for one dimensional Fourier series and were generalized to ultraspherical and Jacobi series by Kogbeliantz and Bavinck in 1925  and 1972, respectively (see also \cite{Wang2000}).
In 1993, Berens and Li \cite{Berens1993} established the relation between the means and the best spherical polynomial approximation on the sphere, and discussed their approximation behavior by various of smoothness. Especially, they proved (see also \cite{Wang2000}) the relation:
\begin{equation}\label{ch1eq6}
\max_{k \geq n} \|V_{k}f-f\|_p\approx \omega \Big(f,\frac1{\sqrt n}\Big)_p, \quad f\in L^p(\sigma).
\end{equation}

Motivated by \cite{Belinsky} and \cite{Ditzian1994}, we will improve the above result. In deed, we will prove
$$
\|V_{n}f-f\|_p \approx \omega \Big(f,\frac1{\sqrt n}\Big)_p
$$
for any $f\in L^p(\sigma)$, $1\leq p\leq +\infty$.

\section{The kennel of de la Vall\'{e}e Poussin}

In the definition of de la Vall\'{e}e-Poussin kernel $v_{n}$ given by (\ref{ch1eq4}), the constants $I_{n,d}$ is requested to satisfy
$$
\int_{\sigma}v_{n}(\stackrel\smile {\mu\nu})d\sigma(\nu)=|\sigma^{d-2}|,
$$
which implies that
$$
\int_{0}^{\pi}v_{n}(\theta)\sin^{2\lambda}\theta d \theta=1 \quad (2\lambda=d-2).
$$

By computation, we have
$$
I_{n,d}=2^{2\lambda}\frac{\Gamma(\lambda+1/2)\Gamma(n+\lambda+1/2)}{\Gamma(n+2\lambda+1)},
$$
where $\Gamma(\lambda)$ is Gamma function. So,
$$
v_n(t)=\frac{\Gamma(n+2\lambda+1)}{2^{2\lambda}\Gamma(\lambda+1/2)\Gamma(n+\lambda+1/2)} \left(\cos\frac {t}2\right)^{2n}.
$$

Since $v_{n}(t)$ are even trigonometric polynomials with degree $n$, $V_{n}(f,\mu)$ are spherical polynomials with degree $n$. So we also call (\ref{ch1eq5}) spherical de la Vall\'{e}e Poussin polynomial operators.

We can translate de la Vall\'{e}e Poussin means given by (\ref{ch1eq5}) into the multiplier form:
\begin{equation}\label{ch2eq1}
V_n(f;\mu)=\sum^{\infty}_{k=0}\omega_{n,k}^{(\lambda)}Y_k(f;\mu)
\end{equation}
where
$$
\omega_{n,k}^{(\lambda)}:=\left\{
\begin{array}{ll}
\displaystyle\frac{n!(n+2\lambda)!}{(n-k)!(n+k+2\lambda)!}, \  \ & 0\leq k\leq n;\\
0, &  k>n.
\end{array}
\right.
$$
Noticing that the means can be rewritten as
\begin{eqnarray*}
V_n(f;\mu)&=&\int^\pi_0 v_n(\theta)S_\theta(f;\mu)\sin^{2\lambda}\theta d\theta\\
&=& \int^\pi_0v_n(\theta) \left(\sum^\infty_{k=0}\frac{P^{\lambda}_k(\cos\theta)}{P^{\lambda}_k(1)} Y_k(f;\mu)\right)\sin^{2\lambda}\theta
d\theta\\
&=&
\sum^\infty_{k=0}\left(\int^\pi_0v_n(\theta)\frac{P^{\lambda}_k(\cos\theta)}{P^{\lambda}_k(1)} \sin^{2\lambda}\theta
d\theta\right)Y_k(f;\mu),
\end{eqnarray*}
it is sufficient to prove
$$
\int_{0}^{\pi}v_{n}(\theta)\frac{P_{k}^{\lambda}(\cos\theta)}{P_{k}^{\lambda}(1)} \sin^{2\lambda}\theta~d\theta = \omega_{n,k}^{(\lambda)}, \quad k \ge 0.
$$
In deed, when $k=1$, one has
\begin{eqnarray*}
&&\int^\pi_0v_n(\theta)\frac{P^{\lambda}_1(\cos\theta)}{P^{\lambda}_1(1)}
\sin^{2\lambda}\theta d\theta\\
&=&\int^\pi_0v_n(\theta)\left(1-2\sin^2\frac{\theta}2\right)
\sin^{2\lambda}\theta d\theta\\
&=&\frac{2^{2\lambda}}{I_{n,d}}\left(\int^\pi_0\left(\cos\frac{\theta}2\right)^{2(n+\lambda)}
\sin^{2\lambda}\frac{\theta}2d\theta-2\int^\pi_0\left(\cos\frac{\theta}2\right)^{2(n+\lambda)}
\sin^{2(\lambda+1)}\frac{\theta}2d\theta\right)\\
&=&\frac{2^{2\lambda+1}}{I_{n,d}}\left(\frac12B\left(\lambda+\frac12,n+\lambda+\frac12\right)-
B\left(\lambda+1+\frac12,n+\lambda+\frac12\right)\right)\\
&=&\frac{n}{n+2\lambda+1}=\frac{n!(n+2\lambda)!}{(n-1)!(n+1+2\lambda)!}=\omega^{(\lambda)}_{n,1},
\end{eqnarray*}
where $B(a,b)$ is Beta function.

Now, we suppose for $k\leq n$,
$$
\int^\pi_0v_n(\theta)\frac{P^{\lambda}_k(\cos\theta)}{P^{\lambda}_k(1)} \sin^{2\lambda}\theta d\theta =\omega^{(\lambda)}_{n,k}.
$$
Then for $k+1$ we first recall the relation (see page 81 of \cite{Szego2003})
$$
(k+1)P^\lambda_{k+1}(x)-2(\lambda+k)xP^\lambda_k(x)+(2\lambda+k-1)P^\lambda_{k-1}(x)=0 \quad (k\ge 1),
$$
i.e.,
$$
P^\lambda_{k+1}(\cos\theta)=\frac1{k+1}\left(2(\lambda+k)\cos\theta P^\lambda_k(\cos\theta) -(2\lambda+k-1)P^\lambda_{k-1}(\cos\theta)\right).
$$
Then,
\begin{eqnarray*}
&&\int^\pi_0v_n(\theta)\frac{P^{\lambda}_{k+1}(\cos\theta)}{P^{\lambda}_{k+1}(1)}
\sin^{2\lambda}\theta d\theta\\
&=& \frac1{P^{\lambda}_{k+1}(1)(k+1)}\left(2(\lambda+k) \int^\pi_0v_n(\theta)\cos\theta P^\lambda_k(\cos\theta)\sin^{2\lambda}\theta d\theta\right.\\
&& \left. -(2\lambda+k-1)\int^\pi_0v_n(\theta)P^\lambda_{k-1}(\cos\theta)\sin^{2\lambda}\theta d\theta\right)\\
&:=& \frac1{P^{\lambda}_{k+1}(1)(k+1)}\left(2(\lambda+k)J_2-J_1\right).
\end{eqnarray*}
By the assumption, we obtain
\begin{eqnarray*}
J_1
&=& (2\lambda+k-1)P^\lambda_{k-1}(1) \int^\pi_0v_n(\theta)\frac{P^\lambda_{k-1}(\cos\theta)}{P^\lambda_{k-1}(1)}\sin^{2\lambda}\theta d\theta\\
&=& \frac{(2\lambda+k-1)!n!(n+2\lambda)!}{\Gamma(2\lambda)(k-1)!(n-k+1)!(n+k-1+2\lambda)!}.
\end{eqnarray*}
For $J_2$ we have
\begin{eqnarray*}
J_2
&=&
\frac{1}{I_{n,d}} \int^\pi_0\left(\cos\frac\theta2\right)^{2n}\left(2\cos^2\frac\theta2-1\right)P^{\lambda}_{k}(\cos\theta)
\sin^{2\lambda}\theta d\theta\\
&=&
\frac{2I_{n+1,d}}{I_{n,d}}P^{\lambda}_{k}(1)\int^\pi_0 v_{n+1}(\theta)
\frac{P^{\lambda}_{k}(\cos\theta)}{P^{\lambda}_k(1)}\sin^{2\lambda}\theta d\theta\\
& &
-P^{\lambda}_k(1)\int^\pi_0 v_{n}(\theta)
\frac{P^{\lambda}_{k}(\cos\theta)}{P^{\lambda}_k(1)}\sin^{2\lambda}\theta d\theta\\
&:=& J_{21}-J_{22},
\end{eqnarray*}
which implies from the assumption that
$$
J_{22}=\frac{\Gamma(k+2\lambda)}{k!\Gamma(2\lambda)}\frac{n!(n+2\lambda)!}{(n-k)!(n+k+2\lambda)!},
$$
and
$$
J_{21}=\frac{2\Gamma(n+1+\lambda+1/2)\Gamma(n+2\lambda+1)}{\Gamma(n+\lambda+1/2)\Gamma(n+1+2\lambda+1)} \frac{\Gamma(k+2\lambda)}{k!\Gamma(2\lambda)}\frac{(n+1)!(n+1+2\lambda)!}{(n+1-k)!(n+1+k+2\lambda)!}.
$$
Therefore,
$$
J_{2} = \frac{\Gamma(k+2\lambda)}{k!\Gamma(2\lambda)} \frac{n!(n+2\lambda)!}{(n+1-k)!(n+1+k+2\lambda)!}\left(n(n+1)+k(2\lambda+k)\right).
$$
So,
\begin{eqnarray*}
\int^\pi_0v_n(\theta)\frac{P^{\lambda}_{k+1}(\cos\theta)}{P^{\lambda}_{k+1}(1)}
\sin^{2\lambda}\theta d\theta
=\frac{n!(n+2\lambda)!}{(n-k-1)!(n+k+1+2\lambda)!}=\omega^{(\lambda)}_{n,k+1}.
\end{eqnarray*}

On the other hand, it is clear that for $k>n$, $\omega^{(\lambda)}_{n,k}=0$. Hence, de la Vall\'{e}e Poussin means $V_n(f;\mu)$ have the form of multiplier expression given in (\ref{ch2eq1}).

Now we give some properties for the de la Vall\'{e}e Poussin kernel $v_n$.

\textbf{Lemma~2.1}~~Let $v_{n}(t)$ be the kernel of de la Vall\'{e}e Poussin defined by (\ref{ch1eq4}), $2\lambda=d-2$ and $d \geq 3$. Then there hold
\begin{equation}\label{ch2eq2}
\int_{0}^{\pi}\theta^{-\lambda}v_n(\theta)\sin^{2\lambda}\theta d \theta \leq C(d)n^{\frac{\lambda}{2}},
\end{equation}
and
\begin{equation}\label{ch2eq3}
\int_{0}^{\pi}\theta^{-\frac{2}{m}}v_n(\theta)\sin^{2\lambda}\theta d \theta \leq C(d)n^{\frac{1}{m}}, \quad m=1,2,\dots.
\end{equation}

\textbf{Proof.}~~We only prove (\ref{ch2eq2}). The proof of (\ref{ch2eq3}) is similar. First, a direct computation implies
$$
I_{n,d}=C(d)\frac{(2n+d-3)!!}{(2n+2d-4)!!}\approx n^{-\frac{d-1}{2}}.
$$
Then,
\begin{eqnarray*}
\int_{0}^{\pi}\theta^{-\lambda}v_n(\theta)\sin^{2\lambda}\theta d\theta
&=& \frac{1}{I_{n,d}}\int_{0}^{\pi}\theta^{-\lambda} \Big(\cos\frac{\theta}{2}\Big)^{2n}\sin^{2\lambda}d\theta\\
&=& \frac{J_{n,d}^{(-\lambda)}}{I_{n,d}},
\end{eqnarray*}
where
$$
J_{n,d}^{(-\lambda)}=\int_{0}^{\frac{\pi}{2}}\theta^{-\frac{d-2}{2}}\sin^{d-2} \theta\cos^{2n}\frac{\theta}{2}d\theta.
$$
So, we have
\begin{eqnarray*}
J_{n,d}^{(-\lambda)}
  &\leq& 2^{\frac{d}{2}}\int_{0}^{\frac{\pi}{2}}\sin^{\frac{d-2}{2}}t\cos^{2n+d-2}tdt\\
  &=& 2^{\frac{d}{2}-1}B\Big(\frac{\frac{d-2}{2}+1}{2},\frac{2n+d-2+1}{2}\Big)=2^{\frac{d}{2}-1}\Gamma(\frac{d}{4})
        \frac{\Gamma(n+\frac{d-1}{2})}{\Gamma(n+\frac{3d-2}{4})}\\
  &=& C(d)\frac{\Gamma(n+\frac{d-1}{2})}{\Gamma(n+\frac{3d-2}{4})} \approx n^{-\frac{d}{4}}.
\end{eqnarray*}
Therefore
$$
\int_{0}^{\pi}\theta^{-\lambda}v_{n}(\theta)\sin^{2\lambda}\theta d \theta
= \frac{J_{n,d}^{(-2)}}{I_{n,d}} \leq C(d) \frac{n^{-\frac{d}{4}}}{n^{-\frac{d-1}{2}}}=C(d)n^{\frac{d-2}{4}}.
   $$
The proof of Lemma~2.1 is completed.\quad$\Box$

\textbf{Lemma~2.2}~~For the kernel of de la Vall\'{e}e Poussin $v_{n}(t)$ defined by (\ref{ch1eq4}),
we have
$$
\int_{0}^{\pi} \theta^{4}v_{n}(\theta)\sin^{2\lambda}\theta d \theta \leq C(d)n^{-2}.
$$

\textbf{Proof.} Since
\begin{eqnarray*}
J_{n,d}^{(4)}
 &=& \int_{0}^{\pi}\theta^{4}\cos^{2n}\frac{\theta}{2} \sin^{2\lambda}\theta d \theta = 2^{d-1}\pi^{4}\int_{0}^{\frac{\pi}{2}}\sin^{d+2}\theta\cos^{2n+d-2}\theta d \theta\\
 &=& \left\{\begin{array}{ll}
       2^{d-2}\pi^{5}\frac{(2n+d-3)!!(d+1)!!}{(2n+2d)!!},  \   \ & \mbox{ if } d
       \mbox{ is even};\\
    2^{d-1}\pi^{4}\frac{(2n+d-3)!!(d+1)!!}{(2n+2d)!!},    \    \  & \mbox { if }  d \mbox{ is odd}
 \end{array}\right.\\
 &=&C(d)\frac{(2n+d-3)!!}{(2n+2d)!!},
\end{eqnarray*}
we have
$$
\int_{0}^{\pi}\theta^{4}v_{n}(\theta)\sin^{2\lambda}\theta d \theta=\frac{J_{n,d}^{(4)}}{I_{n,d}}=C(d) \frac{(2n+2d-4)!!}{(2n+2d)!!} \leq C(d)n^{-2}.
$$
This finishes the proof of Lemma~2.2.\quad$\Box$

\section{Lower bound of approximation for de la Vall\'{e}e Poussin means}

In this section we prove the main result of this paper, which can be stated as follows.

\textbf{Theorem~3.1}~~Let $V_n(f;\mu)$ be  de la Vall\'ee Poussin means on the sphere given by (\ref{ch1eq5}). Then for $f\in L^p(\sigma), 1\leq p\leq +\infty$, there exists a constant $C$ which is independent of $f$ and $n$, such that
$$
\omega\left(f,\frac1{\sqrt n}\right)_p\leq C\|V_nf-f\|_p.
$$

In order to prove the result, we first prove the following lemma.

\textbf{Lemma~3.1}~~For any  $g, Dg, D^2g \in L^{p}(\sigma)$, $1 \leq p \leq \infty$, there exist the constants $A, B$ and $C_2$ which are are independent of $n$ and $g$, such that
$$
\|V_{n}g-g-\alpha(n)Dg\|_p \leq C_2n^{-2}\|D^{2}g\|_p,
$$
where
$0<A/n \leq \alpha(n) \leq B/n$.

\textbf{Proof.}
Since (see (3.6) of \cite{Pawelke1972})
$$
  S_{\theta}(g;\mu)-g(\mu)=\int_{0}^{\theta}\sin^{-2\lambda}t~dt~\int_{0}^{t}\sin^{2\lambda}u~S_{u}(Dg;\mu) du
$$
we have
$$
S_{u}(Dg;\mu)-Dg(\mu)=\int_{0}^{u}\sin^{-2\lambda}\gamma d\gamma \int_{0}^{\gamma}\sin^{2\lambda}\nu ~S_{\nu}(D^{2}g;\mu)d\nu.
$$
Noting that
\begin{eqnarray*}
V_{n}(g;\mu)-g(\mu)
&=& \int_{0}^{\pi}v_{n}(\theta)\sin^{2\lambda}\theta d \theta \int_{0}^{\theta}\sin^{-2\lambda}tdt\int_{0}^{t}
   \sin^{2\lambda}u~S_{u}(Dg; \mu)du\\
&=& Dg(\mu)\int_{0}^{\pi}v_{n}(\theta)\sin^{2\lambda}\theta d\theta~\int_{0}^{\theta}\sin^{-2\lambda}tdt\int_{0}^{t}\sin^{2\lambda}udu\\
& & +\int_{0}^{\pi}v_{n}(\theta)\sin^{2\lambda}\theta d\theta\int_{0}^{\theta}\sin^{-2\lambda}tdt\\
& & \times \int_{0}^{t}\sin^{2\lambda}u \Big(S_{u}(Dg;\mu)-Dg(\mu)\Big)du\\
&:=& Dg(\mu)\alpha(n)+\Psi(g;\mu),
\end{eqnarray*}
where $\alpha(n)=C(d)n^{-1}$ satisfies $0< An^{-1} \leq C(d)n^{-1} \leq Bn^{-1}$, we obtain that from the H\"older-Minkowski's inequality and the contractility of translation operator
\begin{eqnarray*}
\|\Psi g\|_p
&\leq& \|D^{2}g\|_p\int_{0}^{\pi}v_{n}(\theta)\sin^{2\lambda}\theta d\theta\int_{0}^{\theta}\sin^{-2\lambda}tdt
                       \int_{0}^{t}\sin^{2\lambda}udu\\
&    & \times \int_{0}^{u}\sin^{-2\lambda}\gamma d\gamma\int_{0}^{\gamma}\sin^{2\lambda}\nu d\nu\\
& \leq& C_3\|D^{2}g\|_p\int_{0}^{\pi}v_{n}(\theta)\theta^{4}\sin^{2\lambda}\theta d\theta.
\end{eqnarray*}
Thus, from Lemma~2.2 it follows that
$$
\|\Psi g\|_p \leq C_4n^{-2}\|D^{2}g\|_p.
$$
The Lemma~3.1 has been proved.\quad$\Box$

Now we turn to the proof of Theorem~3.1. We first introduce an operator $V^m_n$ given by
$$
V_{n}^{m}(f;\mu)=\sum_{k=0}^{n}\Big(\int_{0}^{\pi}v_{n}(\theta) Q_{k}^{\lambda}(\cos\theta)\sin^{2\lambda}\theta d\theta\Big)^{m} Y_{k}(f;\mu).
$$
Then, form the orthogonality of projection operator $Y_{k}$, it follows that
\begin{eqnarray*}
V_{n}^{m+l}f
&=& \sum_{k=0}^{n}\Big(\int_{0}^{\pi}v_{n}(\theta)Q_{k}^{\lambda}(\cos\theta)
\sin^{2\lambda}\theta d\theta\Big)^{m}\\
& & \times Y_{k}\Big(\sum_{s=0}^{n}\Big(\int_{0}^{\pi}v_{n}(\theta)Q_{s}^{\lambda}(\cos\theta)\sin^{2\lambda}\theta d\theta\Big)^{l}Y_{s}f\Big)\\
&=& V_{n}^{m}(V_{n}^{l}f).
\end{eqnarray*}
Thus, we take $g=V_{n}^{m}f$ and obtain that
\begin{eqnarray*}
\|f-g\|_p = \|f-V_{n}^{m}f\|_p
\leq \sum_{k=1}^{m}\|V_{n}^{k-1}f-V_{n}^{k}f\|_p
\leq m\|f-V_{n}f\|_p,
\end{eqnarray*}
where $V_{n}^{0}f=f.$

Next, we prove the estimation:
$$
  \|DV_{n}^{m}f\|_p \leq \frac{A}{2C_2}C_{1}n\|f\|_p,
$$
where $A$ and $C_2$ are the same as that in Lemma~3.1. In fact, we have
\begin{eqnarray*}
\|DV_{n}^{m}f\|_p
\leq \left\|\sum_{k=0}^{n}k(k+d-2)\Big(\int_{0}^{\pi}v_{n}(\theta)\Big|Q_{k}^{\lambda}(\cos\theta)\Big| \sin^{2\lambda}\theta d\theta \Big)^{m}Y_{k}(f)\right\|_p.
\end{eqnarray*}
Since (see \cite{Belinsky})
$$
|Q_{k}^{\lambda}(\cos\theta)|
\equiv \Big|\frac{P_{k}^{\lambda}(\cos\theta)}{P_{k}^{\lambda}(1)}\Big|
\leq C_{5}\min \Big((k\theta)^{-\lambda},1\Big),
$$
we use (\ref{ch2eq2}) and obtain for $k\theta \geq 1$ and $\theta \leq \pi/2$, that
\begin{eqnarray*}
\|DV_{n}^{m}f\|_p
&\leq& C_{6} \left\|\sum_{k=0}^{n}k(k+d-2)k^{-\frac{d-2}{2}m}\Big(\int_{0}^{\pi}v_{n}(\theta)\theta^{-\lambda}\sin^{2\lambda}\theta d\theta \Big)^{m}Y_{k}(f)\right\|_p\\
&\leq& C_{7}\: n^{\frac{d-2}{4}m}\|f\|_p\sum_{k=0}^{\infty}k^{2-\frac{d-2}{2}m}.
\end{eqnarray*}
For $2-(d-2)m/2 < -1$, i.e. $m > 6/(d-2)$, it is clear that the series $\sum_{k=0}^{\infty}k^{2-\frac{d-2}{2}m}$ is convergence. Thus
$$
\left\|DV_{n}^{m}f\right\|_p \leq C_{8}n^{\frac{d-2}{4}m}\|f\|_p.
$$
For $k\theta \leq 1$, then (\ref{ch2eq3}) implies that
\begin{eqnarray*}
\|DV_{n}^{m}f\|_p
  &\leq& \left\|\sum_{k=0}^{n}\Big(\int_{0}^{\pi}v_{n}(\theta)\theta^{-\frac{2}{m}}(\theta^{2}k(k+d-2))^{\frac{1}{m}}
     \Big|Q_{k}^{\lambda}(\cos\theta)\Big|\sin^{2\lambda}\theta d\theta\Big)^{m} Y_{k}(f)\right\|_p\\
  &\leq& C_{9}\left\|\sum_{k=0}^{n}\Big(\int_{0}^{\pi}v_{n}(\theta)\theta^{-\frac{2}{m}}
     \sin^{2\lambda}\theta d\theta\Big)^{m} Y_{k}(f)\right\|_p\\
  &\leq& C_{10} n \left\|\sum_{k=0}^{\infty}Y_{k}(f)\right\|_p=\frac{A}{2C_2}C_{1}n\|f\|_p,
\end{eqnarray*}
where $A$ and $C_2$ are the same as that in Lemma~3.1. Therefore, when $m > 6/(d-2)$, there holds
$$
\|DV_{n}^{m}f\|_p \leq \frac{A}{2C_2}C_{1}n\|f\|_p.
$$

Without loss generality, we assume $m_{1} > 6/(d-2)$, and $m > 6/(d-2)+m_{1}$ in the rest of the paper.
According to Lemma~3.1 we see that
\begin{eqnarray*}
\alpha(n)\|DV_{n}^{m}f\|_p
&\leq &  \|V_{n}^{m}f-f\|_p+C_2n^{-2}\|D^{2}V_{n}^{m}f\|_p\\
&\leq &  m\|V_{n}f-f\|_p+\frac{AC_{1}}{2}n^{-1}\|DV_{n}^{m-m_{1}}f\|_p\\
&\leq &  m\|V_{n}f-f\|_p+\frac{AC_{1}}{2}n^{-1}\|DV_{n}^{m}f\|_p\\
&     &  +\frac{AC_{1}}{2} n^{-1} \|DV_{n}^{m-m_{1}}(V_{n}^{m_{1}}f-f)\|_p\\
&\leq &  m\|V_{n}f-f\|_p+\frac{AC_{1}}{2n}\|DV_{n}^{m}f\|_p
         +\frac{AC_{1}C_{11}}{2}  \|V_{n}^{m_{1}}f-f\|_p\\
& = &  C_{12}\|V_{n}f-f\|_p+\frac{AC_{1}}{2n}\|DV_{n}^{m}f\|_p.
\end{eqnarray*}
Setting
$\alpha(n)=AC_{1}/n$,
one has
$$
\frac{1}{n}\|DV_{n}^{m}f\|_p\leq \frac{2C_{12}}{AC_{1}}\|V_{n}f-f\|_p.
$$
So from the definition of K-functional it follows
\begin{eqnarray*}
K\left(f,\frac{1}{\sqrt n}\right)
&\leq & \|f-V_{n}^{m}f\|_p+\Big(\frac{1}{\sqrt{n}}\Big)^{2}\|DV_{n}^{m}f\|_p\\
&\leq & m\|f-V_{n}f\|_p+\frac{2C_{12}}{AC_{1}}\|f-V_{n}f\|_p \leq C_{14}\|f-V_{n}f\|_p,
\end{eqnarray*}
which together with (\ref{ch1eq3}) implies
$$
\omega\left(f,\frac{1}{\sqrt n}\right)_p\leq C\|f-V_{n}f\|_p.
$$

This finishes the proof of Theorem~3.1.\quad$\Box$

From (\ref{ch1eq6}) and Theorem~3.1, the following Corollary~3.1 follows directly.

\textbf{Corollary~3.1}~~For any  $f\in L^p(\sigma), 1\leq p\leq \infty$, there holds
\begin{eqnarray*}
\|V_{n}f-f\|_p \approx \omega\left(f,\frac1{\sqrt n}\right)_p.
\end{eqnarray*}

\end{document}